\documentclass[10pt, a4paper]{article}
 \usepackage[cp850]{inputenc}
 \def\be{\begin{equation}}
 \def\ee{\end{equation}}

\usepackage{amssymb}
\usepackage{amsfonts}
\usepackage{graphicx}
\usepackage{amsmath}

\def\<{\langle}
\def\>{\rangle }\def\Tr{\textrm{Tr}\,}
 \def\tr{\textrm{tr}\,}

 =-0.5 in   \headheight =0 cm

 \title{Why Jordan algebras are natural in statistics: quadratic regression implies Wishart distributions}
 \author{ G.
Letac\thanks{Laboratoire de
 Statistique et Probabilit\'es, Universit\'e Paul Sabatier, 31062 Toulouse,
France, e-mail: \texttt{letac@cict.fr}}, J.
Weso{\l}owski\thanks{Wydzia{\l} Matematyki i Nauk Informacyjnych, Politechnika Warszawska, Warszawa, Poland, e-mail: \texttt{wesolo@mini.pw.edu.pl}}}
 \date{\today}

\begin{document}
 \maketitle
 \begin{abstract}If the space $\mathcal{Q}$ of quadratic forms in $\mathbb{R}^n$ is splitted in a direct sum $\mathcal{Q}_1\oplus \ldots\oplus \mathcal{Q}_k$  and if $X$ and $Y$ are independent random variables of $\mathbb{R}^n$, assume that there exist a real number $a$ such that $E(X|X+Y)=a(X+Y)$ and real distinct
numbers $b_1,...,b_k$ such that $E(q(X)|X+Y)=b_iq(X+Y)$ for any  $q$ in $\mathcal{Q}_i.$ We prove that this happens only when $k=2$, when $\mathbb{R}^n$ can be structured in a Euclidean Jordan algebra and when $X$ and $Y$ have Wishart distributions corresponding to this structure. \end{abstract}

\section{Introduction} Let $S_r$ be the set of  $(r,r)$ real symmetric matrices and let $X$ and $Y$ be independent random variables valued in $S_r$ such that they are Wishart distributed $\gamma_{p,\sigma}$ and $\gamma_{p',\sigma}$ , which means that \be\label{lw}\mathbb{E}(e^{-\tr \theta X} )=\det (I_r+\theta \sigma)^{-p}\ee where $\theta$ and $\sigma$  are in the set $P_r$ of the positive definite elements of $S_r$ and $p$ is in \begin{equation}\label{GY1}\Lambda=\{\frac{1}{2},\ldots,\frac{r-1}{2}\}\cup(\frac{r-1}{2},\infty)\end{equation}
(In (\ref{lw}) $\tr$ means trace). Note that for $a=p/(p+p')$ 
\begin{equation}\label{Rao}\mathbb{E}(X|X+Y)=a(X+Y).\ee 
Assume furthermore that $p+p'>\frac{r-1}{2}.$ This implies that $(X+Y)^{-1}$ exists. Then it is known that $Z=(X+Y)^{-1/2}X
(X+Y)^{-1/2}$ and $X+Y$ are independent and that $Z\sim uZu^T$ for any orthogonal $(r,r)$ matrix $u.$ There are many consequences, nuances and characterizations of the Wishart distributions related to this result. One of these consequences is the following fact: for any $s\in S_r$ consider the two quadratic forms on $S_r$  defined by \be \label{q1q2}q_1^s(x)=\frac{1}{2}\tr^2(xs)+\tr( sxsx),\ q_2^s(x)=\tr^2(xs)-\tr (sxsx)\ee and the two numbers 
$$b_1=\frac{p}{p+p'}\ \frac{p+1}{p+p'+1},\ b_2=\frac{p}{p+p'}\ \frac{p-\frac{1}{2}}{p+p'-\frac{1}{2}}.$$ Then for $i=1,2$ and for any $s$ 
\begin{equation}\label{LM}\mathbb{E}(q_i^s(X)|X+Y)=b_iq_i^s(X+Y)\end{equation} This is the particular case $d=1$ of Corollary 2.3 of Letac and Massam (1998). 
 An important fact about this set $(q_1^s,q_2^s)_{s\in S_r}$   is that it spans the whole space of quadratic forms  $\mathcal{Q}$ on $S_r$ (since if $q^s(x)=\tr^2(xs)$ then $\{q^s\ ;\ s\in S_r\}$ spans $\mathcal{Q}$). More specifically denote by $\mathcal{Q}_i$ the subspace of $\mathcal{Q}$ generated by $\{q_i^s\ ;\  s\in S_r\}.$ Then $\mathcal{Q}=\mathcal{Q}_1\oplus\mathcal{Q}_2$ (see for instance Theorem 5.2 below for a proof).

The aim of the paper is to prove a reciprocal statement of (\ref{Rao}) and (\ref{LM}): Let $V$ be a linear real finite dimensional space  (instead of $S_r$) and denote by  $\mathcal{Q}$ the space of all quadratic forms on  $V$.  Fix a decomposition   $\mathcal{Q}=\mathcal{Q}_1\oplus\mathcal{Q}_2\oplus\ldots \oplus\mathcal{Q}_k$ with $k\geq 2$ as a direct sum of linear subspaces. Consider  two independent random variables $X$ and $Y$ with exponential moments satisfying (\ref{Rao}) for some $a$ and $\mathbb{E}(q(X)|X+Y)=b_iq(X+Y)$ for all $q\in \mathcal{Q}_i$ and  for some  distinct real numbers $b_1,\ldots,b_k.$ We show that under these circumstances, necessarily $k=2$ and $X$ and $Y$ are Wishart distributed in the following sense: there necessarily exists a structure of Euclidean Jordan algebra on $V$ (like symmetric matrices, Hermitian matrices, or space with a Lorentz cone) such that $X$ and $Y$ are Wishart on the symmetric cone associated to it. Section 5 contains more detailed information about the two spaces $\mathcal{Q}_1$ and $\mathcal{Q}_2$ of quadratic forms on $S_r$ (or more generally, on a Euclidean Jordan algebra)

\section{Some history of the subject} \textsc{Wishart distributions on $S_r$.} Wishart distributions have been introduced by J. Wishart (1928) as distributions of $Z_1Z_1^T+\cdots+Z_NZ_N^T\sim \gamma_{N/2,2\Sigma}$   where $Z_1,\ldots,Z_N$ are iid in $\mathbb{R}^r$ such that $Z_i\sim N(0,\Sigma).$ Elegant calculations about them are in Bartlett (1933) and the classical reference is Muirhead (1982). For the space  $S_r$ of $(r,r)$ real symmetric matrices the extension of the definition of $\gamma_{p,\sigma}$ from a half integer $p$ to the whole set $\Lambda$ defined by (\ref{GY1}) is made in the fundamental paper of Olkin and Rubin (1962). Proving that a distribution $\gamma_{p,\sigma}$ on the  semi positive definite matrices such that (\ref{lw}) holds only if $p$ is in $\Lambda$ was considered as a challenge by statisticians (see Eaton (1983)) although the appendix of Olkin and Rubin contains an unnoticed proof of it (and unfortunately erroneous: see Casalis and Letac (1994)). This  conjecture was independently proved by Shanbhag (1988) and Peddada and Richards (1989) by quite different means, although a solution already appeared in Gyndikin (1975) and seems to have been well known by analysts, who also call  the set $\Lambda$ and its extensions the Wallach set (see Lassalle (1987) for proofs and references).

\vspace{4mm}
\noindent\textsc{Lukacs-Olkin-Rubin Theorem.} Wishart distributions on $S_r$ are the most natural generalization of the gamma distributions on the positive line. Lukacs (1956) shows  that if $X$ and $Y$ are positive, independent non Dirac random variables and if $Z=X/(X+Y)$, then $Z$ and $X+Y$ are independent if and only if there exists $\sigma,p,p'>0$ such that $X\sim \gamma_{p,\sigma}$ and $Y\sim \gamma_{p',\sigma}.$ This was extended to $S_r$ by Olkin and Rubin (1962) by a proper definition of $Z$ such that $Z$ is symmetric (for instance by choosing $Z=(X+Y)^{-1/2}X
(X+Y)^{-1/2}$ or by choosing $Z=C^{-1}X(C^{-1})^T$ where $C$ is the triangular matrix with positive diagonal elements coming from the Cholesky decomposition $CC^T=X+Y).$ They show that if $X$ and $Y$ are independent non Dirac random semi positive definite matrices in $S_r$ such that $X+Y$ is invertible and such that $Z\sim uZu^T$ for any orthogonal $(r,r)$ matrix $u$ then $Z$ and $X+Y$ are independent if and only if there exists a positive definite matrix $\sigma$ and $p$ and $p'$ in $\Lambda$ with $p+p'>(r-1)/2$ such that $X\sim \gamma_{p,\sigma}$ and $Y\sim \gamma_{p',\sigma}.$ If $Z$ is defined as $(X+Y)^{-1/2}X
(X+Y)^{-1/2}$, Bobecka and Weso{\l}owski (2002) have shown that the invariance hypothesis for $Z$ by the orthogonal group can be dropped provided one assumes that $X$ and $Y$ have smooth densities. Removing this assumption of density is still a challenge. 

\vspace{4mm}
\noindent\textsc{Wishart distributions on Hermitian matrices and on  Euclidean Jordan algebras.} Since normal distributions on Hermitian spaces have been considered  (see e.g. Goodman (1963)), therefore Wishart distributions on Hermitian matrices occur naturally. Actually physicists considered  them quite early (see Mehta (2004)). Carter (1975) in an unpublished PhD thesis extends Olkin and Rubin to this case.

On the other hand, works on the classification of natural exponential families by their variance function have led to the observation that the exponential family $\{\gamma_{p,\sigma};  \sigma\in P_r\}$ of Wishart distributions on $S_r$ with fixed shape parameter $p\in \Lambda$ has a variance function which is  the map from $S_r$ into itself $x\mapsto V(m)(x)=\frac{1}{p}mxm$ where $m$ is in $P_r.$ In other terms, this means that if $\kappa$ is a cumulant function of 
$\gamma_{p,\sigma}$ then for all $x$ in $S_r$ we have $$\kappa''(\theta)(x)=\frac{1}{p}\kappa'(\theta)x\kappa'(\theta).$$ Facts about multivariate distributions such that their corresponding variance functions are quadratic in the mean are collected in Letac (1989). In particular, Wishart distributions obtained from simple Euclidean Jordan algebras are described there. An indispensable  reference for simple Euclidean Jordan algebras is  Faraut and Koranyi (1994) always abreviated F.-K. below. Recall that simple Euclidean Jordan algebras are basically in one to one correspondence with the irreducible symmetric cones (self dual cones in Euclidean space such that the group of automorphisms of the cone acts transitively on it), in the way that $S_r$ is linked to $P_r.$ A quick definition  of the Wishart distribution $\gamma_{p,\sigma}$ on the Jordan algebra $V$ with rank $r$, Peirce constant $d$, cone $\overline{\Omega}$ of square elements, trace and determinant function $\tr$ and $\det$ can be done by its Laplace transform 
$$\int_{\overline{\Omega}}e^{-\tr \theta x}\gamma_{p,\sigma}(dx)=\det(e+\theta \sigma)^{-p}$$
where $\sigma$ is in the interior $\Omega$ of $\overline{\Omega}$ and where $p$ is in the Gyndikin set of the Jordan algebra $V$ defined by \be \label{GY2}\Lambda_V=\{\frac{d}{2},d, \ldots,\frac{d}{2}(r-1)\} \cup (\frac{d}{2}(r-1),\infty).\ee While the definition of determinant is the standard one for $S_r$ and for Hermitian matrices, it requires some care for the three other types of Jordan algebras: quaternionic Hermitian matrices, 27 dimensional Albert algebra and the algebra of the Lorentz cone.

 Particular cases of use of Wishart distributions on Jordan algebras in statistics occurred earlier (Andersson (1975) for the Hermitian and quaternionic cases, and Jensen (1988) for the Lorentz cone, with its deep connexions to Clifford algebras). Jordan algebras are the natural framework for Wishart  distributions: Casalis and Letac (1996) is  a clarification and an extension to Jordan algebras of Olkin and Rubin (1962) and of  Carter (1975); Carter follows step by step the  difficult Olkin and Rubin's approach and his work was unknown to  Casalis and Letac (1996).  

\vspace{4mm}
\noindent\textsc{Quadratic homogeneity and Wishart distributions.}
A remarkable fact about the classical Wishart distributions on $S_r$ is that the above variance function $m\mapsto V(m)$ is not only quadratic in $m$ but homogeneous quadratic. This happens also  to be true for Wishart distributions on any  Euclidean Jordan algebra. This observation lead Casalis (1991) to prove the converse: any natural   exponential family with a homogeneous quadratic variance function  is a Wishart family, as conjectured in Letac (1989). Put in other words, if $\kappa$ is a cumulant function of some random variable $X$ valued in $\mathbb{R}^n$  such that $\kappa''(\theta)=V(\kappa'(\theta))$ where $V$ is a homogeneous quadratic function, then $\mathbb{R}^n$ can be structured in a Jordan algebra such that $X$ is Wishart for that structure.

\vspace{4mm}
\noindent\textsc{Quadratic regression property.} A slight extension of Lukacs (1956) is to take two non Dirac independent rv $X$ and $Y$ on the positive line such that there exist positive $a$ and $b$ such that $\mathbb{E}(X|X+Y)=a(X+Y)$ and $\mathbb{E}(X^2|X+Y)=b(X+Y)^2$ and to prove that there exist positive $p,p',\sigma$ such that $X\sim \gamma_{p,\sigma}$ and $Y\sim \gamma_{p',\sigma}.$ To see this, just multiply these two equalities by $e^{\theta(X+Y)},$
take expectations and obtain two differential equations for the Laplace transforms of $X$ and $Y.$ This procedure is contained  in Laha and Lukacs (1960). Bivariate regression version of Lukacs theorem based on
conditions $E(X_i^2|X+Y)=b(X_i+Y_i)^2,$ $i=1,2,$ where $X=(X_1,X_2)$ and
$Y=(Y_1,Y_2)$ are independent was obtained in Wang (1981). This result was
generalized in Letac and Wesolowski (2008) by considering regressions of
quadratic forms
$E(q(X)|X+Y)=b q(X+Y)$ for all quadratic forms $q$ orthogonal to an arbitrary
fixed quadratic form $q_0.$ That is in the setting of the present paper we
required $k=1$ and codimension of $\mathcal{Q}_1$ to be equal 1.

Letac and Massam (1998) use the quadratic regression  approach to get a simpler proof of Olkin and Rubin theorem, as extended to Jordan algebras in Casalis and Letac (1996). It actually characterizes the Wishart distributions of independent $X$ and $Y$ in $S_r$ (and more generally of a Jordan algebra) through the following  properties: if for $i=1,2,$ $s\in S_r$ and  $q^s_i$ are defined by (\ref{q1q2}),  then (\ref{LM}) holds
(with suitable analogues of $q_i$ if the Jordan algebra is not $S_r).$ Note that this regression perspective leads to a characterization of $\gamma_{p,\sigma},\  \gamma_{p',\sigma}$ without the hypothesis  of invertibility of $X+Y$ which was needed in the Olkin and Rubin characterization. 

\section{Main result} Let $V$ be a real linear space with dimension $n>1$, let $V^*$ be its dual and consider the space $\mathcal{F}=L_s(V,V^*)$ of the symmetric linear maps from $V$ to $V^*$. If $\theta\in V^*$ and $x\in V$ we write $\<\theta, x\>$ for $\theta(x).$ Denote by $\mathcal{Q}$ the space of quadratic forms $q$ on $V,$ namely the set of real functions $q$ on $V$ such that $(x,y)\mapsto \frac{1}{2}(q(x+y)-q(x)-q(y))$ is bilinear on $V\times V$ and $q(\lambda x)=\lambda^2q(x)$ when $\lambda$ is a real number. The map from $\mathcal{F}$ to $\mathcal{Q}$  defined by  $f\mapsto q_f$ where  $x\mapsto q_f(x)=\<f(x),x\>$ is one to one. More specifically:
$$\frac{1}{2}(q_f(x+y)-q_f(x)-q_f(y))=\frac{1}{2}(\<f(x),y\>+\<f(y),x\>)=\<f(x),y\>$$ 
For $q\in \mathcal{Q}$ we therefore define the inverse map $q\mapsto f_q$ of $f\mapsto q_f$ by 
$$\frac{1}{2}(q(x+y)-q(x)-q(y))=\<f_q(x),y\>.$$ 

Let us also define here the concept of irreducibility for a probability measure $\mu$ on $V.$  We say first that $\mu$ is reducible if there exists a  direct sum $V_1\oplus V_2=V$ with $\dim V_i>0$ for $i=1,2$, two probability measures $\mu_1$ and $\mu_2$ on $V_1$ and $V_2$ such that $\mu=\mu_1\otimes \mu_2.$ In other terms, if $X\sim \mu$ its projections $X_1$  on $V_1$ parallel to $V_2$ and $X_2$ on $V_2$ parallel to $V_1$ are independent.   Suppose that furthermore $X$ has a Laplace transform $L=e^{\kappa}$ defined on some open set $\Theta\subset V^*=V_1^*\oplus V_2^*.$ In this case $\kappa(\theta)=\kappa_1(\theta_1)+ \kappa_2(\theta_2)$ where $\theta_i$ is the projection of $\theta$ on $V_i^*$ and $\kappa_1$ and $\kappa_2$ are the cumulant functions of $X_1$ and $X_2.$ We say also that $X$ and $\kappa$ are reducible in that case. Finally, $\mu$, $X$ and $\kappa$  are said to be irreducible if they are not reducible...

\vspace{4mm}\noindent \textbf{Theorem 3.1} 
Let $\mathcal{Q}_1\oplus \mathcal{Q}_2\oplus\ldots \oplus \mathcal{Q}_k=\mathcal{Q}$ be a direct sum decomposition of the space of  quadratic forms on $V$ with $k\geq 2.$ Let $X$ and $Y$ be two independent irreducible random variables valued in $V$ such that their Laplace transforms exist on an open set $\Theta\subset V^*.$ We assume that \begin{enumerate}\item there exists a real number $a$ such that $\mathbb{E}(X|X+Y)=a(X+Y);$ 
\item there exist distinct numbers $b_1,\ldots, b_k$  such that for any $i=1,\ldots ,k$ and for any $q\in \mathcal{Q}_i$ we have 
\be \label{qr}\mathbb{E}(q(X)|X+Y)=b_iq(X+Y).\ee
\end{enumerate}
Under these circumstances $0<a<1$, $k=2$ and  there exists a simple Euclidean Jordan algebra structure on $V$ such that $X$ and $Y$ are Wishart distributed on the positive cone of the algebra with the same scale parameter and respective shape parameters $p$ and $p'$  in   $\Lambda_V$  defined in (\ref{GY2}). Moreover $ \mathcal{Q}_1$ and  $ \mathcal{Q}_2$ are spanned by \be \label{qs} q_1^s(x)=\frac{d}{2}\tr^2(xs)+\tr( \mathbb{P}(x)(s)s),\ q_2^s(x)=\tr^2(xs)-\tr( \mathbb{P}(x)(s)s)\ee where $\tr$, $\mathbb{P}$ and $d$ are respectively the trace, the quadratic map and the Peirce constant of the Jordan algebra and $s\in V.$ In this case  
\be\label{abb}a=\frac{p}{p+p'},\ b_1=\frac{p}{p+p'}\ \frac{p+1}{p+p'+1},\ b_2=\frac{p}{p+p'}\ \frac{p-\frac{d}{2}}{p+p'-\frac{d}{2}}.\ee

\vspace{4mm}\noindent \textbf{Proof.} 
Denote by $L_X$ and $L_Y$ the Laplace transforms of $X$ and $Y.$ It is standard to prove that from condition 1) we have  $L_X^{1-a}=L_{Y}^{a}$: just multiply both sides of  $\mathbb{E}(X|X+Y)=a(X+Y)$ by $e^{\<\theta,X+Y\>}$ where $\theta\in \Theta$ and take expectations of both sides to obtain the differential equation $aL'_X/L_X=(1-a)L'_Y/L_Y.$ The fact that $X$ and $Y$ are irreducible implies that $a=0$ or $a=1$ is impossible. The fact that $\log L_X$ and $\log L_Y$ are  convex implies that $a<0$ or $ a>1$ are impossible.  From now on  we denote $e^{\kappa}=L_X=L_{Y}^{a/(1-a)}.$ 

In the sequel, we use the symbol $\Tr$ for the trace of an endomorphism. The symbol $\tr$ is reserved for the trace in a Jordan algebra. If $q$ is a quadratic form on $V$ we write 
$$q(\frac{\partial}{\partial \theta})(\kappa)(\theta)=\Tr (f_q \kappa''(\theta)).$$  Since $\kappa$ is a real twice differentiable function defined on an open subset of $V^*,$ the second derivative $\kappa''(\theta)$ is an element of $L_s(V^*,V)$, the linear map $f_q$ is an element of $L_s(V,V^*)$ and thus $f_q \kappa''(\theta)$ belongs to $ L(V^*,V^*)$.  It therefore makes sense to speak of the trace of  this endomorphism of $V^*.$ Note that $\<f_q(x),x\>=\Tr (f_q (x\otimes x))$ and that $\frac{\partial}{\partial \theta}\otimes \frac{\partial}{\partial \theta} \kappa=\kappa''.$  This explains the definition $q(\frac{\partial}{\partial \theta})(\kappa)=\Tr (f_q \kappa'').$ Also $q(\kappa')$ can be  written in terms of $f_q$ as $q(\kappa')=\Tr (f_q (\kappa'\otimes \kappa'))=\<f_q(\kappa'),\kappa'\>.$

Calculations done in Letac and Weso{\l}owski (2008) (2.9), show that for any $i=1,\ldots k$ and for all $q\in \mathcal{Q}_i$ we have 
\begin{equation}\label{llo}(1-\frac{b_i}{a})q(\frac{\partial}{\partial \theta})(\kappa)=(\frac{b_i}{a^2}-1)q(\kappa').\end{equation} (Again, to prove  (\ref{llo}) just multiply (\ref{qr}) by $e^{\<\theta,X+Y\>}
$ and take expectations). Observe that $b_i=a$ is impossible, since it implies that $q(\kappa')=0$ for any $q$ in $\mathcal{Q}_i.$ Since $\mathcal{Q}_i$ is not the zero space, there exists a non zero $q$ with $q(\kappa')=0.$ Now $\{x\in V; q(x)=0\}$ is a quadric of $V$ and has an empty interior. On the other hand, since $X$ is irreducible, this implies that $X$ cannot be concentrated on some affine subspace of $V$. Therefore $\kappa$ is strictly convex and the set $\kappa'(\Theta)$ is open and cannot be contained in a quadric. Thus $a=b_i$ is impossible, division by $(1-\frac{b_i}{a})$ is permitted and we rewrite (\ref{llo}) as 
\begin{equation}\label{yky}q(\frac{\partial}{\partial \theta})(\kappa)=p_iq(\kappa')\end{equation} where $p_i=\frac{b_i-a^2}{a^2-ab_i}.$ 

Now let us fix  $\theta\in V^*$  and consider the element $\theta\otimes \theta$ of $\mathcal{F}$ defined by $(\theta\otimes \theta) (x)=\<\theta, x\> \theta.$ Denote by $\mathcal{F}_i$ the image of $\mathcal{Q}_i$ by the isomorphism $q\mapsto f_q.$ Obviously we have $$\mathcal{F}_1\oplus \mathcal{F}_2\oplus\ldots \oplus \mathcal{F}_k=\mathcal{F}.$$ Therefore there exist elements $f_i\in \mathcal{F}_i$ such that $f_1+\ldots+f_k=\theta\otimes \theta.$ Since $f_1,\ldots,f_k$ depend actually on $\theta$ we rather write $f_i(\theta,x)$ instead of $f_i(x)$ for $x\in V.$ Thus $x\mapsto f_i(\theta,x)$ is a linear map from $V$ to $V^*$. 
We rewrite the equality $\theta\otimes \theta=f_1+\ldots+f_k$ as 
$$\<\theta,x\>^2=\<f_1(\theta,x),x\>+\ldots+\<f_k(\theta,x),x\>$$ for any $x$ in $V.$ We now fix $\theta=\theta_0$ in this equality and we recall that $q(\frac{\partial}{\partial \theta})(\kappa)(\theta)$ means $\Tr (f_q \kappa''(\theta)).$  Thus we apply the equality $$\<\theta_0,\frac{\partial}{\partial \theta}\>^2=\<f_1(\theta_0,\frac{\partial}{\partial \theta}),\frac{\partial}{\partial \theta}\>+\ldots+\<f_k(\theta_0,\frac{\partial}{\partial \theta}),\frac{\partial}{\partial \theta}\>$$
to $\kappa$, the log of $L.$ We get
$$\Tr((\theta_0\otimes \theta_0)\ \kappa''(\theta))=\sum_{i=1}^k\Tr[f_i(\theta_0,\cdot)\kappa''(\theta)].$$
We now use the fact that $x\mapsto \<f_i(\theta_0,x),x\>=q(x)$ is a quadratic form belonging to $\mathcal{Q}_i$ to which we apply (\ref{yky}). We therefore get 
\begin{equation}\label{ykz}\Tr((\theta_0\otimes \theta_0)\  \kappa''(\theta))=\sum_{i=1}^k p_i\Tr(f_i(\theta_0,\cdot)\kappa'(\theta)\otimes \kappa'(\theta)).\end{equation}
Since this is true for any $\theta_0$ in $V^*$ this is enough to claim that $\kappa''$ is a quadratic homogeneous function of $\kappa'.$

We now apply the  Casalis' theorem (1991) which says that if $\kappa$ is irreducible and if $\kappa''$ is a quadratic homogeneous  function of $\kappa'$, then there exists a simple Euclidean Jordan algebra structure on $V$ related to $X$ in a way that we explain now.  Let $\Omega$ be the open cone of the squares of $V$, let $\tr$ and $\det$ be the trace and determinant functions on the Jordan algebra, let $d$ and $r$ be the Peirce and rank constants of $V.$  
Then there exists $p\in \Lambda_V$ defined by (\ref{GY2}) and $\sigma \in \Omega$ such $X$ has the Wishart distribution $\gamma_{p,\sigma}$ on $\overline{\Omega}$
defined by its Laplace transform  $\mathbb{E}(e^{-\tr \theta X} )=\det (I_r+\theta \sigma)^{-p}$ for all $\theta\in \Omega.$ 
 
To complete the proof, denote for a while  by $ \widetilde{\mathcal{Q}}_1$ and  $\widetilde{ \mathcal{Q}}_2$ the spaces of quadratic forms spanned by $(q_1^s)_{s\in V}$ and $(q_2^s)_{s\in V}$ as defined in (\ref{qs}). Denote also 
$$\widetilde{b}_1=\frac{p}{p+p'}\ \frac{p+1}{p+p'+1},\ \widetilde{b}_2=\frac{p}{p+p'}\ \frac{p-\frac{d}{2}}{p+p'-\frac{d}{2}}.$$

Recall that we want to prove that $k=2$ and that $\{\widetilde{\mathcal{Q}}_1,\widetilde{\mathcal{Q}}_2\}=\{\mathcal{Q}_1,\mathcal{Q}_2\}.$
Let now $q\in \mathcal{Q}_i$. Therefore $\mathbb{E}(q(X)|X+Y)=b_iq(X+Y).$ We now write $q=q_1+q_2$ with $q_i\in \widetilde{\mathcal{Q}}_i$ which is possible since $\widetilde{\mathcal{Q}}_1\oplus \widetilde{\mathcal{Q}}_2=\mathcal{Q}.$ Recall that since $X$ and $Y$ have distributions $\gamma_{p,\sigma}$ and $\gamma_{p',\sigma}$  we can write $\mathbb{E}(q_i(X)|X+Y)=\widetilde{b}_iq(X+Y).$ Thus 
$$(\widetilde{b}_1-b_i)q_1(X+Y)=(b_i-\widetilde{b}_2)q_2(X+Y).$$
Since $X+Y$ is valued in the open set $\Omega$ this implies $(\widetilde{b}_1-b_i)q_1=(b_i-\widetilde{b}_2)q_2.$
Thus the two sides of this equality are zero: either $b_i=\widetilde{b}_1$ and $q_2=0$ or the reverse statement holds. Since we have assumed that $b_1,\ldots,b_k$ are distinct, this ends the proof. 
\section{Comments}\begin{enumerate}\item Surprisingly enough, while starting from a linear space $V$ without any additional algebraic structure, the regression conditions on $X$ and $Y$ of the theorem impose by themselves a Euclidean Jordan algebra structure on $V.$ 

\item The three numbers $a$, $b_1$ and $b_2$ together with   the dimension of $V$ determine uniquely the structure of Jordan algebra on $V$ in the following sense: we can see from the equations (\ref{abb}) that $b_2<a^2<b_1<a.$ Moreover these equations  give the Peirce constant $d$ of $V$ by 
$$d=2\frac{a-b_1}{b_1-a^2}\, \frac{a^2-b_2}{a-b_2}.$$
Since the rank $r$ satisfies $\dim V=r+\frac{d}{2}r(r-1)$ the type of the Jordan algebra is completely known. 

\item In the theorem, $k=1$ would lead to $X$ and $Y$ concentrated on a line $\mathbb{R}v$ of $V$. If $X=X_1v$ and $Y=Y_1v$ then $X_1$ and $Y_1$ would be one dimensional gamma distributed and $X$ would not be irreduciblesince we have assumed  $\dim V>1.$ Furthermore if in the theorem we do not assume that $b_1,\ldots,b_k$ are distinct, then either they are all equal to one $b$ and this sends us back to the trivial case $k=1$ or they are not and if $k'\geq 2$ is the number of distinct $b_i$'s, then the theorem gives $k'=2.$ 

\item Some comments about irreducibility are in order. If $L_Y$ is a power of $L_X$, then $Y$ is irreducible if and only if $X$ is. Therefore irreducibility can  be assumed in the theorem  for $X$ only. If irreducibility is not assumed, we have an artificial generality. For instance suppose that $(X_1,X_2,X_3,Y_1,Y_2,Y_3)$ are independent real rv with $X_i\sim \gamma_{\alpha_i,\sigma}$ and $Y_i\sim \gamma_{\beta_i,\sigma}.$ Then for $i\neq j$ we have 
\begin{eqnarray*}\mathbb{E}(X_iX_j|X+Y)&=& \frac{\alpha_i\alpha_j}{(\alpha_i+\beta_i)(\alpha_j+\beta_j)}(X_i+Y_i)(X_j+Y_j),\\  \mathbb{E}(X_i^2|X+Y)&=& \frac{\alpha_i(\alpha_i+1)}{(\alpha_i+\beta_i)(\alpha_i+\beta_i+1)}(X_i+Y_i)^2.\end{eqnarray*}
This implies that $k=6$ corresponding to the 6 independent quadratic forms on $V=\mathbb{R}^3$ defined by  $q_{ij}(x)=x_ix_j$ for $i\leq j.$ 
\end{enumerate}
\section{The spaces $\mathcal{Q}_1$ and $\mathcal{Q}_2$: the operator $\Psi$ }
If $V$ is a simple Euclidean Jordan algebra with rank $r$ and Peirce constant $d=2d'$,  denote by $\mathcal{F}=L_s(V)$ the space of symmetric linear operators on $V.$ The dimension of $V$ is $n=r+dr(r-1)/2$. 
Given $y\in V$, important examples of elements of $\mathcal{F}$ are respectively  $\mathbb{L}(y)$  defined by $x\mapsto xy$ where $xy$ is the Jordan product, 
and $$\mathbb{P}(y)=2(\mathbb{L}(y))^2-\mathbb{L}(y^2)$$ as defined  in F.-K. page 32. If $a$ and $b$ are in $V$ we denote by $a\otimes b$  the  endomorphism $x\mapsto a\, \tr(bx)$ of $V.$ The endomorphism  $a\otimes b+b\otimes a$ belongs to $\mathcal{F}.$ 
Denote by $\mathcal{F}_1$ and $\mathcal{F}_2$ the linear subspaces of $\mathcal{F}$  respectively generated by  and $\{d'y\otimes y+\mathbb{P}(y); y\in V\}$ and $\{y\otimes y-\mathbb{P}(y); y\in V\}.$  From (\ref{qs}) $\mathcal{F}_1$ and $\mathcal{F}_2$ are canonically isomorphic to $\mathcal{Q}_1$ and $\mathcal{Q}_2$ by $q\mapsto f_q$ where $q(x)=\<f_q(x),x\>.$ We endow $\mathcal{F}$
with the Euclidean structure defined by $\Tr (ab).$ Here again we distinguish the trace $\tr$ of the Jordan algebra $V$ from the trace $\Tr$ of the endomorphisms on the linear space $V.$ Here is a list of various traces:

\vspace{4mm}\noindent \textbf{Proposition 5.1}\begin{enumerate} 
\item $\Tr(a\otimes b)=\tr(ab),$ 
$\Tr [(a\otimes b)(c\otimes d)]=\tr(ad)\tr(bc)$, 
 $$\Tr((a_1\otimes b_1)\cdots (a_k\otimes b_k))=\tr(a_1b_k)\tr(a_2b_1)\cdots \tr(a_kb_{k-1}).$$
\item $\Tr[\mathbb{L}(a)\mathbb{L}(b)(c\otimes d)]=\tr [(a(bc))d]$

\item 
$\Tr(\mathbb{P}(a)(b\otimes c))=\tr[(\mathbb{P}(a)b) c]$

\end{enumerate}
\vspace{4mm}\noindent \textbf{Proof} 1) is standard since it only involves the Euclidean structure of $V$ and not its Jordan algebra structure. 2) is a consequence of 1). Applying the definition of  $\mathbb{P}(a)$, 3) is a consequence of 2).

In the theorem below, we consider an endomorphism $\Psi$ of $\mathcal{F}$ such that $\Psi(y\otimes y)=\mathbb{P}(y)$ for all $y\in V.$ It is an essential tool of the two papers Casalis and Letac (1996) and Letac and Massam (1998). The theorem shows that $\mathcal{F}_1$ and $\mathcal{F}_2$ are its two eigenspaces and uses this fact to give the dimensions of the  spaces of quadratic forms $\mathcal{Q}_1$ and $\mathcal{Q}_2 $ defined
in Th. 3.1 above.

\vspace{4mm}\noindent \textbf{Theorem 5.2}
\begin{enumerate} \item There exists a symmetric  endomorphism $\Psi$ of $\mathcal{F}$ such that $\Psi(y\otimes y)=\mathbb{P}(y)$ for all $y\in V.$ It satisfies  \be\label{clm}\Psi (\mathbb{P}(y))=d'y\otimes y+(1-d')\mathbb{P}(y)\ee\item The spaces $\mathcal{F}_1$ and $\mathcal{F}_2$ are orthogonal and  $\mathcal{F}=\mathcal{F}_1\oplus \mathcal{F}_2$   \item The spaces $\mathcal{F}_1$ and $\mathcal{F}_2$ are the two eigenspaces of $\Psi$  corresponding to the two eigenvalues  $1$ and $-d'$ respectively. \item The dimensions of $\mathcal{F}_1$ and $\mathcal{F}_2$ are given by $$\frac{n(n+1)}{2}-\dim \mathcal{F}_1=\dim \mathcal{F}_2=\frac{r(r-1)}{2}\times \frac{1+d'(2r-3)+d'^2(r-1)(r-2)}{1+d'}$$ \end{enumerate}
\vspace{4mm}\noindent \textbf{Examples.}
For the Jordan algebra  associated to the Lorentz cone where $r=2$ we get $\dim \mathcal{F}_2=1.$ More specifically, if $E$ is a Euclidean space with scalar product $\vec{x}.\vec{y}$ consider the Jordan algebra $V=\mathbb{R}\times E$ endowed with the Jordan product between $x=(x_0,\vec{x})$ and $y=(y_0,\vec{y})$ 
defined by $$xy=(x_0y_0+\vec{x}.\vec{y}\, ,\,  x_0\vec{y}+y_0\vec{x}).$$  Here the Lorentz cone is $ \{(x_0,\vec{x})\in V\ ;\  x_0>\|\vec{x}\|\},$ the trace is  $\tr (x_0,\vec{x})=2x_0$ and the Peirce constant is $d=\dim E -1.$ In this case $\mathcal{F}_2$ is spanned by the symmetry $S$ defined by $(x_0,\vec{x})\mapsto (x_0,-\vec{x}).$ To see this observe that if $e=(1,\vec{0})$ then $S=e\otimes e-\mathbb{P}(e)$ is in $\mathcal{F}_2$ and use $\dim \mathcal{F}_2=1.$ As a consequence if $\left[\begin{array}{cc}a&b\\b^*&c\end{array}\right]$ represents a symmetric endomorphism of $V$ (where $a$ is real, $c$ is a symmetric endomorphism of $E$ and $b$ is a linear form on $E$) then $\left[\begin{array}{cc}a&b\\b^*&c\end{array}\right]$  is in 
$\mathcal{F}_1$ if and only if it is orthogonal to $$S=\left[\begin{array}{cc}1&0\\0&-\mathrm{id}_E\end{array}\right],$$ that is if and only if $a=\Tr c.$ The dimension of $\mathcal{F}_1$ is $\frac{1}{2}(n-1)(n+2).$

For the Jordan algebra $S_r$ of symmetric real matrices where $d=1,$ we get $\dim \mathcal{F}_2=\frac{r^2}{12}(r-1)(r+1)$ and  $\dim \mathcal{F}_1=\frac{r}{24}(r+1)(r^2+5r+6).$ For the Jordan algebra of  Hermitian  matrices where $d=2,$ we get $$\dim \mathcal{F}_1=\left(\frac{r(r+1)}{2}\right)^2,\ \dim \mathcal{F}_2=\left(\frac{r(r-1)}{2}\right)^2,$$  and since $d'=1,$ $\Psi$ is an orthogonal symmetry with respect to $\mathcal{F}_2.$ For the Jordan algebra of Hermitian quaternionic  matrices where $d=4,$ we get $\dim \mathcal{F}_2=4r\frac{r(r-1)(r-2)}{6}+\frac{r(r-1)}{2}$ and $\dim \mathcal{F}_1=\frac{r^2}{3}(4r^2-1).$  For the Albert algebra where $d=8$ and $r=3$ we get $\dim \mathcal{F}_2=27,\ \dim \mathcal{F}_1=351=27\times 13.$

\vspace{4mm}\noindent \textbf{Proof.} 1) The existence of $\Psi$ is proved in Casalis and Letac (1996) (Lemma 6.3) and (\ref{clm}) is proved in Letac and Massam (1998) (Proposition 3.1). For proving  that $\Psi$ is symmetric, enough is to see that $\Tr[\Psi(x\otimes x)(y\otimes y)]$ is symmetric in $x$ and $y$ in $V$ since $\{y\otimes y\ ; \ y\in V\}$ spans $\mathcal{F}.$ Equivalently we have to see that $\Tr[\mathbb{P}(x)(y\otimes y)]$ is symmetric. From Proposition 5.1 part 3, we have to show that  
$\tr[(\mathbb{P}(x)y) y)]$ is symmetric. Applying the definition of $\mathbb{P},$ we get
$$\tr[(\mathbb{P}(x)y) y)]=\tr[(2(x(xy)-x^2y)y].$$ Let us now use Proposition II.1.1, (iii) in
 F.-K. which says $$\mathbb{L}(x^2y)-\mathbb{L}(x^2)\mathbb{L}(y)=2\mathbb{L}(xy)\mathbb{L}(x)-2\mathbb{L}(x)\mathbb{L}(y)\mathbb{L}(x).$$ Applying this equality to $y$ we get $(x^2y)y-x^2y^2=2(xy)^2-2x(y(xy))$ that we rewrite as $2(xy)^2+x^2y^2=2x(y(xy))+(x^2y)y.$ Since the left hand side is symmetric in $(x,y)$ this proves $2x(y(xy))+(x^2y)y=2y(x(xy))+(y^2x)x$ which implies in turn that $(2(x(xy)-x^2y)y$ is symmetric in $(x,y)$ and shows that  $\Psi$ is symmetric. 

2) and 3) Since $\{y\otimes y\ ;\  y\in V\}$ spans $\mathcal{F}$ and since 
$$y\otimes y=\frac{1}{1+d'}(d'y\otimes y+\mathbb{P}(y))+\frac{1}{1+d'}(y\otimes y-\mathbb{P}(y)),$$ clearly $\mathcal{F}=\mathcal{F}_1+ \mathcal{F}_2.$ From the formula (\ref{clm}) and the definition of $\Psi$
 we get easily that $\mathcal{F}_1$ and $ \mathcal{F}_2$ are made of eigenvectors of $\Psi$ respectively for the eigenvalues $1$ and $-d'.$ In particular $\mathcal{F}_1\cap \mathcal{F}_2=\{0\}.$ Therefore ${\mathcal F}={\mathcal F}_1\oplus {\mathcal F}_2$ and thus the endomorphism $\Psi$ has no other eigenvalues. From the fact that $\Psi$ is symmetric, $\mathcal{F}_1$ and $\mathcal{F}_2$ are orthogonal.

4) It is the difficult point. We have  $\dim \mathcal{F}_1+\dim \mathcal{F}_2=\frac{n(n+1)}{2}$ where $n$ is the dimension of $V.$ An other linear equation for $(\dim \mathcal{F}_1,\ \dim \mathcal{F}_2)$ is 
$\mathrm{trace} (\Psi)= \dim \mathcal{F}_1-d'\dim \mathcal{F}_2$ leading to \be\label{TR}\dim \mathcal{F}_2=\frac{1}{1+d'}\left(\frac{n(n+1)}{2}-\mathrm{trace} (\Psi)\right).\ee We embark for a calculation of $\mathrm{trace} (\Psi)$ by selecting an orthonormal basis $f=(f_{\ell})_{\ell=1}^{n(n+1)/2}$ of $\mathcal{F}$ and by computing $\Tr[\Psi(f_{\ell})f_{\ell}]$ in order to get 
$$\mathrm{trace} (\Psi)=\sum_{\ell=1}^{n(n+1)/2}\Tr[\Psi(f_{\ell})f_{\ell}].$$
The basis $f$ is chosen as follows. We start from a Jordan frame $(c_1,\ldots,c_r)$ of $V$ (see F.-K. page 44).  Recall
that $c_s^2=c_s$ and $c_s c_t=0$ for $s\ne t.$
We denote by  
$V(c,\lambda)$ the eigenspace of $V$ of  $\mathbb{L}(c)$ for the eigenvalue $\lambda$. For $1\leq s< t\leq r$ we denote 
$$V_{st}=V(c_s,\frac{1}{2})\cap V(c_t,\frac{1}{2}),\ V_{ss}=V(c_s,1).$$
Recall that $V=\bigoplus_{1\leq s\leq t\leq r}V_{st}$, that the dimension of $V_{st}$ is $d$ for $s<t$ and $1$ for $s=t$  and that these spaces are orthogonal (F.-K. Th. IV 2.1, (i)). Let $(c_{s,t}^1,\ldots,c_{s,t}^d)$  be    an orthonormal basis  of the space $V_{s,t}$ for $s<t$. The space $V_{ss}$ is spanned by $c_s.$ For simplicity denote also by $e=(e_1,\ldots,e_n)$ the orthonormal basis of $V$ defined by the $c_s'$s and the $c_{st}^k$'s. 
Finally the basis $f$ of $\mathcal{F}$ consists of the elements of the form $f_{\ell}=e_i\otimes e_i$ for $i=1,\ldots,n$, or $f_{\ell}=(e_i\otimes e_j+e_j\otimes e_i)/\sqrt{2}$ for $1\leq i<j\leq n$. 
Since $e$ is an orthonormal basis of the Euclidean space $V$ it is standard to see  that $f$ is an orthonormal basis of the space $\mathcal{F}$ of symmetric endomorphisms of $V.$

We now compute $\Tr[\Psi(f_{\ell})f_{\ell}]=C_{\ell}$ for all possible choices of $f_{\ell}$ in the basis $f.$ 
\begin{enumerate}\item \textsc{Case A:} $f_{\ell}=e_i\otimes e_i.$ From Proposition 5.1, part 5 we have for all $x\in V:$
\be\label{Aa}\Tr(\mathbb{P}(x)(x\otimes x))=\tr x^4\ee

Case A1: $e_i=c_s$. Thus inserting $x=c_s$ in (\ref{Aa}) we get 
$C_{\ell}=\tr c_s=1.$

Case A2: $e_i=c_{st}^k$. We use the fact that $x^2=\frac{\|x\|^2}{2}(c_s+c_t)$ when $x\in V_{st}$ (see F.-K. Proposition IV. 1.4 (i)) and apply (\ref{Aa}) to $x=c_{st}^k.$ We get 
$$C_{\ell}=\tr(c_{st}^k)^4 =\frac{1}{4}\tr[(c_t+c_s)^2]=\frac{1}{2}.$$

\item \textsc{Case B:} $f_{\ell}=(e_i\otimes e_j+e_j\otimes e_i)/\sqrt{2}.$
We use the following calculation: 
$$\Psi(x\otimes y+y\otimes x)=\mathbb{P}(x+y)-\mathbb{P}(x)-\mathbb{P}(y)=2[\mathbb{L}(x)\mathbb{L}(y)+\mathbb{L}(y)\mathbb{L}(x)-\mathbb{L}(xy)]$$
(F.-K. page 32) and, using Proposition 5.1 part 2: 
$$\Tr[(\mathbb{L}(x)\mathbb{L}(y)+\mathbb{L}(y)\mathbb{L}(x)-\mathbb{L}(xy))(x\otimes y+y\otimes x)]=\tr[(yx^2)y+(xy^2)x].$$ Thus \be \label{BC}C_{\ell}=\tr[(e_je_i^2)e_j+(e_ie_j^2)e_i].\ee

Case B1: $e_i=c_s,\ e_j=c_t$ with $s<t.$ From (\ref{BC}):
$$C_{\ell}=\tr [(c_tc_s^2)c_t+(c_sc_t^2)c_s]=0.$$

Case B2: $e_i=c_s,\ e_j=c_{uv}^k$ with $1\leq u<v\le r$ with $s\in\{u,v\}.$ By the definition of $V_{uv}$ we have $c_{uv}^kc_s=\frac{1}{2}c_{uv}^k$ and thus from (\ref{BC}):
\begin{eqnarray*}C_{\ell}&=&\tr [(c_{uv}^kc_s^2)c_{uv}^k+(c_s(c_{uv}^k)^2)c_s]=\tr [(c_{uv}^kc_s)c_{uv}^k+(c_s\frac{1}{2}(c_u+c_v))c_s]\\&=&\tr [\frac{1}{4}(c_u+c_v)+(c_s\frac{1}{2}(c_u+c_v))c_s]=1.\end{eqnarray*}

Case B3: $e_i=c_s,\ e_j=c_{uv}^k$ with $1\leq u<v\le r$ with $s\not \in\{u,v\}.$ Here we have $c_{uv}^kc_s=0$ from F.-K. page 68 last line. A calculation similar to B2 gives $
C_{\ell}=0.$

Case B4: $e_i=c_{uv}^k,\ e_j=c_{uv}^m$ with  $1\le u<v\le r$ and $1\le k<m\le d.$ Note that if $x$ and $y$ have norm 1 in $V_{uv}$ then 
$$(yx^2)y+(xy^2)x=(y\frac{1}{2}(c_u+c_v))y+(x\frac{1}{2}(c_u+c_v))x=\frac{1}{2}(y^2+x^2)=\frac{1}{2}(c_u+c_v)$$
Applying this to $x=c_{uv}^k$ and $y=c_{uv}^m$ we get  $
C_{\ell}=1$ through (\ref{BC}).

Case B5: $e_i=c_{st}^k,\ e_j=c_{uv}^m$ with  $1\le s<t\le r$,  with  $1\le u<v\le r$, with $1\le k,m\le d$ and  with $\{s,t\}\cap\{u,v\}$ reduced to one point, say $u=s.$ Note that if $x$ and $y$ have norm 1 in $x\in V_{st}$ and $y\in V_{sv}$ then 
$$(yx^2)y+(xy^2)x=(y\frac{1}{2}(c_s+c_t))y+(x\frac{1}{2}(c_s+c_v))x=\frac{1}{4}(y^2+x^2)=\frac{1}{8}(2c_s+c_t+c_v)$$
Applying this to $x=c_{st}^k$ and $y=c_{sv}^m$ we get  $
C_{\ell}=1/2$ through (\ref{BC}). 

Case B6: $e_i=c_{st}^k,\ e_j=c_{uv}^m$ with  $1\le s<t\le r$ and  with  $1\le u<v\le r$ with $\{s,t\}\cap\{u,v\}=\emptyset.$ Using  F.-K. page 68 last line we see that $(yx^2)y+(xy^2)x=0$ when $x\in V_{st}$ and $y\in V_{uv}$. Therefore $
C_{\ell}=0.$

\end{enumerate}

We are now in position to compute the trace  of $\Psi.$ We adopt the obvious notation $C(A_1)=\sum_{\ell \in A_1} C_{\ell}.$
Thus $$\mathrm{trace} (\Psi)=C(A_1)+C(A_2)+C(B_2)+C(B_4)+C(B_5).$$ Since $C_{\ell}$ is constant on each of these five sets $A_1,A_2,B_2,B_4,B_5$ we first count the number of their elements: 
$$N(A_1)=r,\ N(A_2)=r(r-1)d',\ N(B_2)=2r(r-1)d',$$$$\ N(B_4)=r(r-1)d'(d'-\frac{1}{2}),\ N(B_5)=2r(r-1)(r-2)d'^2.$$ 
We get finally
$$\mathrm{trace} (\Psi)=r+r(r-1)d'\left[2+(r-1)d'\right]$$ which leads to the result through (\ref{TR}).

\vspace{4mm}\noindent \textbf{Comments.} Observe that $\Tr [\Psi(f)f]=\Tr f^2$ if and only if $f$ is in $\mathcal{F}_1$ (write $f=f_1+f_2$ with $f_i\in \mathcal{F}_i$ and $\Tr [\Psi(f)f]-Tr f^2=-(d'+1)\Tr f_2^2$ to see this). Thus in the above orthonormal basis $(f_{\ell})_{1}^{n(n+1)/2}$ of $\mathcal{F}$ we have $f_{\ell}\in \mathcal{F}_1$ if and only if $C_{\ell}=1$, which happens only in the cases A1, B2 and B4. This is a set of size $N(A_1)+N(B_2)+N(B_4)<\dim \mathcal{F}_1.$ Similarly $\Tr [\Psi(f)f]=-d'\Tr f^2$ if and only if $f$ is in $\mathcal{F}_2,$ and this shows that no $f_{\ell}$ is in $\mathcal{F}_2$ since $C_{\ell}\geq 0$ for all $\ell.$
\section{References}\vspace{4mm}\noindent
\textsc{Andersson, S.} (1975)
 'Invariant normal models.'
 {\it Ann. Statist.} {\bf
3 }, 132-154.

\vspace{4mm}\noindent
\textsc{Bartlett, D. J.} (1933) 'On the theory of the statistical regression.' {\it Proc. 
R. Soc. Edimb.}, {\bf 53}, 260-283.

\vspace{4mm}\noindent \textsc{Bobecka, K.  and Weso{\l}owski, J.} (2002) 'The Lukacs-Olkin-Rubin theorem without the invariance of the ''quotient''.'
 \textit{Studia Mathematica} \textbf{152}, 147-160.

\vspace{4mm}\noindent
\textsc{Carter, E.} (1975) \textit{'Characterization and testing problems in the complex Wishart distribution.' } Ph. D. thesis, University of Toronto.

\vspace{4mm}\noindent 
\textsc{Casalis, M.} (1991)
 'Les familles exponentielles \`a variance quadratique homog\`ene sont des lois de Wishart sur un c\^one sym\'etrique.'
 {\it C. R. Acad. Sci. Paris S\'er. I Math.} {\bf
278}, 293-295.

\vspace{4mm}\noindent
\textsc{Casalis, M. and Letac, G.} (1994)
 'Characterization of the Jorgensen set in the generalized linear model.'
 {\it Test} {\bf
3}, 145-162.

\vspace{4mm}\noindent
\textsc{Casalis, M. and Letac, G.} (1996)
 'The Lukacs-Olkin-Rubin characterization of the Wishart distribution on symmetric cones.'
 {\it Ann. Statist.} {\bf
24 }, 763-786.

\vspace{4mm}\noindent \textsc{Eaton, M.L.} (1983)   
{\it Multivariate Statistics: A Vector-Space Approach.} John Wiley, New York.

\vspace{4mm}\noindent\textsc{Faraut, J. and  Kor\'anyi, A. } (1994) {\it
Analysis on Symmetric Cones.} Oxford University Press.

\vspace{4mm}\noindent \textsc{Gyndikin, S.} (1975)   'Invariant generalized functions in homogeneous spaces.' 
{\it J. Funct. Anal. Appl.,} \textbf{9}, 50-52.

\vspace{4mm}\noindent\textsc{Goodman, N.R.} (1963) `Statistical 
analysis based on a certain multivariate complex Gaussian distribution.' 
{\it Ann. Math. Statist.} {\bf 34} 152-176.

\vspace{4mm}\noindent \textsc{Jensen, S. T.} (1988) 'Covariance hypotheses which are linear in both the 
covariance and the inverse covariance.' {\it Ann. Statist.} {\bf 16,} 302-322.

\vspace{4mm}\noindent
\textsc{Laha, R.G. and Lukacs, E.} (1960)' On a problem connected to quadratic regression.' \textit{Biometrika} \textbf{47}, 335-343.
 
\vspace{4mm}\noindent
\textsc{Lassalle, M.} (1987) 'Alg\`ebre de Jordan et ensemble de Wallach.' \textit{Invent. math.} \textbf{89}, 375-393.

\vspace{4mm}\noindent
\textsc{Letac, G.} (1989) 'Le probl\`eme de la classification des familles exponentielles naturelles
 de $\mathbb{R}^d$ ayant une fonction variance quadratique.'  Probability on Groups IX, Oberwolfach.
 {\it Lecture notes in mathematics,} Springer {\bf
1379}, 192-216.

\vspace{4mm}\noindent
\textsc{Letac, G. and   Massam, H.} (1998)
 'Quadratic and inverse regression for Wishart distributions.'
 {\it Ann.  Statist.} {\bf
26}, 573-595.

\vspace{4mm}\noindent
\textsc{Letac, G.  and Weso{\l}owski, J.} (2008) 'Laplace transforms which are negative powers of quadratic polynomials.' {\it  Trans. Amer. Math. Soc.}, \textbf{360}, 6475-6496. 

\vspace{4mm}\noindent
\textsc{Lukacs, E.} (1955)
 'A characterization of the gamma distribution.'    
 \textit{Ann. Math. Statist.}. {\bf
26}, 319-324.

\vspace{4mm}\noindent
\textsc{Mehta, M. L.} (2004) \textit{'Random matrices'} Pure and Applied Math. Series \textbf{142}, Elsevier, London.

\vspace{4mm}\noindent\textsc{Muirhead, R.J.} (1982) {\it Aspects of
Multivariate Analysis.} Wiley, New York.

\vspace{4mm}\noindent\textsc{Olkin, I. and Rubin, H.} (1962) 'A characterization of the Wishart 
distribution.' {\it Ann. Math. Statist., 33}, 1272-1280.

\vspace{4mm}\noindent\textsc{Peddada S.D. and Richards D. St. P.} (1991)  'Proof of a conjecture of M.L. Eaton on the characteristic function of the Wishart distribution.' {\it Ann. Probab.} \textbf{19}, 869-874. 
272-280.

\vspace{4mm}\noindent\textsc{Shanbhag, D.N. } (1988) 'The Davidson
Kendall problem and related results on the structure of the
Wishart distribution.' {\it Austr. J. Statist.}, {\bf 30A}
272-280.

\vspace{4mm}\noindent\textsc{Wishart, J. } (1928) 'The generalised product moment distribution in samples from a normal multivariate population.' \textit{Biometrika}, \textbf{20A} No 1/2, 32-52.

\vspace{4mm}\noindent \textsc{Wang, Y.} (1981) Extensions of
Lukacs' characterization of the gamma distribution. In: {\it
Analytic Methods in Probability Theory}, Lect. Notes in Math. {\bf
861}, Springer, New York, 166-177.

\end{document}